\input amstex
\documentstyle{amsppt}
\pagewidth{5.4in}
\pageheight{7.6in}
\magnification=1200
\TagsOnRight
\NoRunningHeads
\topmatter
\title
\bf A harmonic map flow associated with\\ 
the standard solution of Ricci flow
\endtitle
\author
Shu-Yu Hsu
\endauthor
\affil
Department of Mathematics\\
National Chung Cheng University\\
168 University Road, Min-Hsiung\\
Chia-Yi 621, Taiwan, R.O.C.\\
e-mail:syhsu\@math.ccu.edu.tw
\endaffil
\date
Feb 7, 2007
\enddate
\address
e-mail address:syhsu\@math.ccu.edu.tw
\endaddress
\abstract
Let $(\Bbb{R}^n,g(t))$, $0\le t\le T$, $n\ge 3$, be a standard solution of 
the Ricci flow with radially symmetric initial data $g_0$. We will extend 
a recent existence result of P.~Lu and G.~Tian and prove that for any 
$t_0\in [0,T)$ there exists a solution of the corresponding harmonic map 
flow $\phi_t:(\Bbb{R}^n,g(t))\to (\Bbb{R}^n,g(t_0))$ satisfying $\partial
\phi_t/\partial t=\Delta_{g(t),g(t_0)}\phi_t$ of the form $\phi_t(r,\theta)
=(\rho (r,t),\theta)$ in polar coordinates in $\Bbb{R}^n\times (t_0,T_0)$, 
$\phi_{t_0}(r,\theta)=(r,\theta)$, where $r=r(t)$ is the radial 
co-ordinate with respect to $g(t)$ and $T_0=\sup\{t_1\in (t_0,T]:
\|\widetilde{\rho}(\cdot ,t)\|_{L^{\infty}(\Bbb{R}^+)}
+\|\partial\widetilde{\rho}/\partial r(\cdot ,t)\|_{L^{\infty}(\Bbb{R}^+)}
<\infty\quad\forall t_0<t\le t_1\}$ with $\widetilde{\rho}(r,t)
=\log (\rho(r,t)/r)$. We will also prove the uniqueness of solution of the 
harmonic map flow within the class of functions of the form $\phi_t(r,\theta)
=(\rho (r,t),\theta)$, $\rho (r,t)=re^{\widetilde{\rho}(r,t)}$, for some 
function $\widetilde{\rho}(r,t)$. We will also use the same technique to 
prove that the
solution $u$ of the heat equation in $(\Omega\setminus\{0\})\times (0,T)$ has
removable singularities at $\{0\}\times (0,T)$, $\Omega\subset\Bbb{R}^m$, 
$m\ge 3$, if and only if $|u(x,t)|=O(|x|^{2-m})$ locally uniformly on
every compact subset of $(0,T)$.
\endabstract
\endtopmatter
\NoBlackBoxes
\define \pd#1#2{\frac{\partial #1}{\partial #2}}
\define \1{\partial}
\define \2{\overline}
\define \3{\varepsilon}
\define \4{\widetilde}
\define \5{\underline}
\document

It is known that Ricci flow is a powerful method in studying the geometry of
manifolds. A manifold $(M, g(t))$, $0\le t\le T$, with an evolving metric 
$g(t)$ is said to be a Ricci flow if it satisfies
$$
\frac{\1 }{\1 t}g_{ij}=-2R_{ij}
$$
in $M\times (0,T)$. Short time existence of solution of Ricci flow on compact 
manifold was proved by R.~Hamilton \cite{H1} using the Nash-Moser Theorem. 
Short time existence of solutions of the Ricci flow on complete
non-compact Riemannian manifold with bounded curvature was proved by W.X~Shi
\cite{S1}. Global existence and uniqueness of solutions of the Ricci flow 
on non-compact manifold $\Bbb{R}^2$ was obtained by S.Y.~Hsu in \cite{Hs1}.
We refer the readers to the lecture notes by B.~Chow \cite{C} and the book
\cite{CK} by B.~Chow and D.~Knopf on the basics of Ricci flow. Interested
readers can also read the papers of R.~Hamilton \cite{H1-6}, S.Y.~Hsu 
\cite{Hs1-7}, B.~Kleiner and J.~Lott \cite{KL}, J.~Morgan and 
G.~Tian \cite{MT}, G.~Perelman \cite{P1}, \cite{P2}, W.X.~Shi 
\cite{S1}, \cite{S2}, L.F.~Wu \cite{W1}, \cite{W2}, R.~Ye \cite{Ye} for some 
of the most recent results on Ricci flow.

Since the proof of existence of solution of the Ricci flow in \cite{H1} is
very hard and there is very few uniqueness results for Ricci flow, later 
D.M.~DeTurck \cite{D} deviced another method to prove existence and uniqueness
of solution of Ricci flow. For any $t_0\in [0,T)$ he introduced an auxillary 
harmonic map flow $\phi_t:(M,g(t))\to (M,g(t_0))$ associated with the Ricci 
flow given by
$$
\frac{\1}{\1 t}\phi_t=\Delta_{g(t),g(t_0)}\phi_t\tag 0.1
$$
where
$$
\Delta_{g(t),g(t_0)}\phi_t=\Delta_{g(t)}\phi_t+g^{ij}(x,t)
\Gamma^{\alpha}_{\beta ,\gamma}(\phi_t(x))\frac{\1\phi_t(x)^\beta}{\1 x^i}
\frac{\1\phi_t(x)^\gamma}{\1 x^j}
$$
in the local co-ordinates $x=(x^1,\dots,x^n)$ of the domain manifold
$(M,g(t))$ and the local co-ordinates $(y^{\alpha})$ of the target manifold 
$(M,g(t_0))$ with $\Gamma^{\alpha}_{\beta ,\gamma}$ being the Christoffel 
symbols of $(M,g(t_0))$. This harmonic map flow then induces a push forward
metric $\hat{g}(t)=(\phi_t)_{\ast}(g(t))$ on the target manifold $M$ which 
satisfies the Ricci-DeTurck flow \cite{H5}
$$
\frac{\1 }{\1 t}\hat{g}_{\alpha\beta}=(L_V\hat{g})_{\alpha\beta}
-2\hat{R}_{\alpha\beta}
$$
for some time varying vector field $V$ on the target manifold $M$ 
where $\hat{R}_{\alpha\beta}$
is the Ricci curvature associated with the metric $\hat{g}(t)$. The existence 
and uniqueness of solutions of Ricci flow on compact manifolds are then 
reduced to the study of existence and other properties of the harmonic map 
flow (0.1) and the Ricci-DeTurck flow \cite{H5}.

In \cite{P1}, \cite{P2}, G.~Perelman proposed a scheme to study Ricci 
flow with singularities. Essential to this scheme is the construction of 
a standard solution of the Ricci flow which is used to replace the solution
near the singularities during surgery. In \cite{LT} P.~Lu and G.~Tian proved
the short time existence of solution of the harmonic map flow associated 
with the standard solution $(\Bbb{R}^n,g(t))$, $0\le t\le T$, $n\ge 3$, 
of Ricci flow with radially symmetric initial data $g_0$. 

In this paper we will extend their result and prove that for any 
$t_0\in [0,T)$ there exists a solution of the corresponding harmonic map 
flow $\phi_t:(\Bbb{R}^n,g(t))\to (\Bbb{R}^n,g(t_0))$ satisfying 
$$\left\{\aligned
&\1\phi_t/\1 t=\Delta_{g(t),g(t_0)}\phi_t\quad\text{ in }
\Bbb{R}^n\times (t_0,T_0)\\
&\phi_{t_0}(x)=x\qquad\qquad\qquad\text{ in }\Bbb{R}^n
\endaligned\right. \tag 0.2
$$ 
of the form 
$$
\phi_t(r,\theta)=(\rho (r,t),\theta),\rho (r,t)=re^{\4{\rho}(r,t)},\tag 0.3
$$ 
for some function $\4{\rho}(r,t)$ in polar coordinates in $\Bbb{R}^n\times 
(t_0,T_0)$, 
$\phi_{t_0}(r,\theta)=(r,\theta)$, where $r=r(t)$ is the radial co-ordinates 
with respect to the metric $g(t)$ and
$$
T_0=\sup\{t_1\in (t_0,T]:\|\widetilde{\rho}(\cdot ,t)\|_{L^{\infty}(\Bbb{R}^+)}
+\|\1\widetilde{\rho}/\1 r(\cdot ,t)\|_{L^{\infty}(\Bbb{R}^+)}
<\infty\quad\forall t_0<t\le t_1\}.\tag 0.4
$$
Then
$$
\4{\rho}(r,t)=\log\biggl (\frac{\rho(r,t)}{r}\biggr ).\tag 0.5
$$
By (0.2),
$$
\rho (r,t_0)=r\quad\forall r>0.\tag 0.6
$$
We will also prove the uniqueness of solution of the harmonic map flow
(0.2) within the class of functions of the form (0.3) for some function 
$\4{\rho}(r,t)$.

The plan of the paper is as follows. In section one we will extend the 
existence result of \cite{LT} and prove the existence of solution of (0.2) 
of the form (0.3) in $\Bbb{R}^n\times (t_0,T_0)$ where $T_0$ is given by 
(0.4). In section two we will prove various estimates for the Green function 
of the heat equation in cylindrical and punctured cylindrical domains. In 
section three we will use the Green function estimates to prove that the 
transformed solution in $(\Bbb{R}^{n+2}\setminus\{0\})\times (t_0,T_0)$ 
has removable singularities on the line $\{0\}\times (t_0,T_0)$. We will 
also prove the uniqueness of solution of (0.2) in section three. In section
four we will prove that a solution $u$ of the heat equation in 
$(\Omega\setminus\{0\})\times (0,T)$ has removable singularities at 
$\{0\}\times (0,T)$, $\Omega\subset\Bbb{R}^m$, $m\ge 3$, if and only if 
there exists $\2{B_R(0)}\subset\Omega$ such that
$$
|u(x,t)|=O(|x|^{2-m})\quad\text{ uniformly on }[t_1,t_2]\quad\forall
0<|x|\le R, 0<t_1<t_2<T.\tag 0.7
$$ 
We first start with a definition. Let $n\ge 3$. For any $0\le t<(n-1)/2$, let 
$h(t)$ be the standard metric on $S^{n-1}$ with constant scalar curvature 
$$
\frac{1}{1-\frac{2t}{n-1}}.
$$ 
Let $g_0$ be a fixed rotationally symmetric complete smooth metric with 
non-negative curvature operator on $\Bbb{R}^n$ such that 
$(\Bbb{R}^n\setminus\2{B(0,2)},g_0)$ is isometric to the half infinite 
cylinder $(S^{n-1}\times\Bbb{R}^+,h(1)\times ds^2)$ (cf. Section 1 of 
\cite{LP} and definition 12.1 of \cite{MT}). By the argument of section 1 
of \cite{LP} such $g_0$ exists.  We say that a Ricci flow $(\Bbb{R}^n,g(t))$, 
$0\le t<T$, is a standard solution if $g(0)=g_0$ and the curvature $Rm$ is 
locally bounded in time $t\in [0,T)$. By the results of \cite{LP}, there 
exists a standard solution $(\Bbb{R}^n,g(t))$ of the Ricci flow on $(0,T)$
for some $T\in (0,\frac{n-1}{2})$ with $g(0)=g_0$ which has non-negative 
curvature operator $Rm(t)$ for each $t\in [0,T)$. 

We will now let $(\Bbb{R}^n,g(t))$, $0\le t<T$, be the standard solution 
of Ricci flow for the rest of the paper. By the result of \cite{P2} and 
\cite{LT} for each $0\le t<T$, $g(t)$ is a rotationally symmetric metric of
$\Bbb{R}^n$. Let $\hat{r}$ be the standard radial co-ordinate in $\Bbb{R}^n$.
As observed by P.~Lu and G.~Tian \cite{LT} if $d\sigma =h(1)$ is 
the standard metric on $S^{n-1}$ with constant scalar curvature 1, then there 
exists a function $f(r,t)\ge 0$ such that
$$
g(t)=dr^2+f(r,t)^2d\sigma
$$
where $r=r(\hat{r},t)$ is the radial co-ordinate on $\Bbb{R}^n$ with respect
to the metric $g(t)$. We fix a $t_0\in [0,T)$ and consider the harmonic map 
flow $\phi_t:(\Bbb{R}^n,g(t))\to (\Bbb{R}^n,g(t_0))$ (0.2) of the form (0.3). 
Let $\4{\rho}(r,t)$ be given by (0.5) and let $f_0(r)=f(r,t_0)$. Then by
(0.5) and (0.6),
$$
\4{\rho}(r,t_0)=0.
$$
By the computation in \cite{LT} and \cite{MT}, 
$$
f(r,t)=re^{\4{f}(r^2,t)}, f_0(\rho)=\rho\, e^{\4{f}_0(\rho^2)}
$$
for some smooth functions $\4{f}(w,t)$ and $\4{f}_0(w)$ of $w\ge 0$ and $t$
with $\4{f}(w,t_0)=\4{f}_0(w)$. Moreover $\4{\rho}(r,t)$ satisfies
$$
\frac{\1\4{\rho}}{\1 t}=\frac{\1^2\4{\rho}}{\1 r^2}+\frac{n+1}{r}
\frac{\1\4{\rho}}{\1 r}+\biggl [(n-1)\frac{\1\4{f}}{\1 r}-r\xi\biggr ]
(r^2,t)\frac{\1\4{\rho}}{\1 r}
+\biggl (\frac{\1\4{\rho}}{\1 r}\biggr )^2+G(\4{\rho},r^2,t)=0\tag 0.8
$$
where
$$\align
G(\4{\rho},w,t)=&\frac{n-1}{w}\biggl [1-e^{2\4{f}_0(\rho^2)-2\4{f}(w,t)}
\biggr ]+2(n-1)\frac{\1\4{f}}{\1 w}(w,t)\\
&\qquad -2(n-1)e^{2\4{f}_0(\rho^2)+2\4{\rho}-2\4{f}(w,t)}
\frac{\1\4{f}_0}{\1 w}(\rho^2)-2\xi (w,t)\endalign
$$
and $\xi (w,t)$ is a smooth function satisfying
$$
\xi (r^2,t)=\frac{1}{r}\frac{\1 r}{\1 t}.
$$
Let $x=(x^1,\dots,x^{n+2})\in\Bbb{R}^{n+2}$, 
$|x|=(\sum_{i=1}^{n+2}(x^i)^2)^{1/2}$, and 
$$
\2{\rho}(x,t)=\4{\rho}(|x|,t),\quad\2{f}(x,t)=\4{f}(|x|^2,t), 
\quad G_1(s, x,t)=G(s,|x|^2,t).\tag 0.9
$$
If $\4{\rho}\in C^{2,1}([0,\infty)\times (t_0,T'))$ is a solution of (0.8) 
in $\Bbb{R}^+\times (t_0,T')$ for some $T'\in (t_0,T)$ \cite{LT}, then 
$\2{\rho}(x,t)$ is a radially symmetric solution of 
$$
\frac{\1\2{\rho}}{\1 t}=\Delta\2{\rho}
+\nabla [(n-1)\2{f}-\2{B}]\cdot\nabla\2{\rho}+|\nabla\2{\rho}|^2
+G_1(\2{\rho},x,t)\tag 0.10
$$
in $(\Bbb{R}^{n+2}\setminus\{0\})\times (t_0,T')$ where
$$
B(w,t)=\frac{1}{2}\int_0^{w}\xi (u,t)\,du
$$
is a smooth even function of $w$ and $\2{B}(x,t)=B(|x|^2,t)$ is a smooth 
function of $x\in\Bbb{R}^{n+2}$. 

Conversely if $\2{\rho}(x,t)$ is a radially 
symmetric solution of (0.10) in $\Bbb{R}^{n+2}\times (t_0,T')$, then 
$\4{\rho}\in C^{2,1}([0,\infty)\times (t_0,T'))$ is a solution of (0.8) 
in $\Bbb{R}^+\times (t_0,T')$.

As in \cite{LT} we rewrite (0.10) as
$$
\frac{\1\2{\rho}}{\1 t}=\Delta\2{\rho}+F(x,\2{\rho},\nabla\2{\rho},t)
\tag 0.11
$$
where
$$
F(x,\2{\rho},\nabla\2{\rho},t)=\nabla [(n-1)\2{f}-\2{B}](x,t)\cdot
\nabla\2{\rho}+|\nabla\2{\rho}|^2+G_1(\2{\rho},x,t).
$$
We will fix $t_0\in (0,T)$ for the rest of the paper.

For any set $A$ we let $\chi_A$ be the characteristic function of the
set $A$. For any $R>0$, $m\ge 2$, let $B_R=\{x:|x|<R\}\subset\Bbb{R}^m$ and 
let $G_R=G_R(x,t,y,s)$ be the Green function of the heat equation in 
$B_R\times (-\infty,\infty)$. Then (\cite{LSU} P.408) 
$$
G_R(x,t,y,s)=\Gamma (x,t,y,s)-g_R(x,t,y,s)\quad\forall x,y\in\2{B}_R, t>s,
$$
where 
$$
\Gamma (x,t,y,s)=\frac{1}{(4\pi (t-s))^{\frac{m}{2}}}
e^{-\frac{|x-y|^2}{4(t-s)}}
$$
and $\forall y\in\2{B}_R, s\in\Bbb{R}$,
$$\left\{\aligned
&\1_tg_R=\Delta_xg_R\qquad\qquad\qquad\forall x\in B_R, t>s\\
&g_R(x,t,y,s)=\Gamma (x,t,y,s)\quad\forall x\in\1 B_R,t>s\\
&g_R(x,s,y,s)=0\qquad\qquad\quad\forall x\in\1 B_R.
\endaligned\right.\tag 0.12
$$
Note that by the maximum principle (cf. \cite{A},\cite{F}),
$$
0\le G_R(x,t,y,s)\le G_{R'}(x,t,y,s)\le\Gamma (x,t,y,s)\quad\forall
x,y\in B_R,0<R<R',t>s.\tag 0.13
$$

$$
\text{Section 1}
$$

In this section we will extend the existence result of P.~Lu and G.~Tian 
\cite{LT} and prove the existence of solution of (0.2) of the
form (0.3) in $\Bbb{R}^n\times (t_0,T_0)$ for some constant $T_0\in (t_0,T]$
given by (0.4). We first recall a result of \cite{LT}.

\proclaim{\bf Theorem 1.1}(Section 2.2.3 of \cite{LT})
There exists $T'\in (t_0,T]$ such that (0.11) has a radially symmetric 
solution $\2{\rho}$ in $\Bbb{R}^{n+2}\times (t_0,T')$ satisfying
$$
\2{\rho}(x,t_0)\equiv 0\quad\text{ in }\Bbb{R}^{n+2},\tag 1.1
$$
$$
\2{\rho}(x,t)=\int_{t_0}^t\int_{\Bbb{R}^{n+2}}
\frac{1}{(4\pi (t-s))^{\frac{n+2}{2}}}e^{-\frac{|x-y|^2}{4(t-s)}}
F(y,\2{\rho},\nabla\2{\rho},s)\,dy\,ds\tag 1.2
$$
for all $x\in\Bbb{R}^{n+2}$, $t_0\le t<T'$, and
$$
\sup_{t_0\le t<T'}(\|\2{\rho}(\cdot,t)\|_{L^{\infty}(\Bbb{R}^{n+2})}
+\|\nabla\2{\rho}(\cdot,t)\|_{L^{\infty}(\Bbb{R}^{n+2})})<\infty.\tag 1.3
$$
\endproclaim

\proclaim{\bf Theorem 1.2}
There exists $T_0\in (t_0,T]$ such that (0.11) has a radially symmetric 
solution $\2{\rho}$ in $\Bbb{R}^{n+2}\times (t_0,T_0)$ satisfying (1.1), 
(1.2), in $\Bbb{R}^{n+2}\times (t_0,T_0)$ where $T_0$ is given by (0.4) 
with $\2{\rho}$ and $\4{\rho}$ being related by (0.9).
\endproclaim
\demo{Proof}
By Theorem 1.1 there exists $T'\in (t_0,T]$ such that (0.11) has a radially 
symmetric solution $\2{\rho}$ in $\Bbb{R}^{n+2}\times (t_0,T')$ satisfying 
(1.1) and (1.2). Let $T_1\in (t_0,T]$ be the maximal existence time of a 
radially symmetric solution $\2{\rho}$ of (0.11) in $\Bbb{R}^{n+2}\times 
(t_0,T)$ that satisfies (1.1), (1.2), and (1.3) for any $t_0<T'<T_1$. 
Then $T_1\le T_0$. We claim that $T_1=T_0$. Suppose $T_1<T_0$. Then
$$
C_1=\sup_{t_0\le t<T_1}(\|\2{\rho}(\cdot,t)\|_{L^{\infty}(\Bbb{R}^{n+2})}
+\|\nabla\2{\rho}(\cdot,t)\|_{L^{\infty}(\Bbb{R}^{n+2})})<\infty.
$$
Let $T_2\in (t_0,T_1)$ be a constant to be determined later.
We will now use a modification of the argument of \cite{LiT} and \cite{LT}
to construct a solution of (0.11) with initial data $\2{\rho}(x,T_2)$.
For any $T_2\le t<T$, $x\in\Bbb{R}^{n+2}$, let $\2{\rho}_1(x,t)
=\2{\rho}(x,T_2)$, and 
$$
\2{\rho}_i(x,t)=v(x,t)+\int_{T_2}^t\int_{\Bbb{R}^{n+2}}
\frac{1}{(4\pi (t-s))^{\frac{n+2}{2}}}e^{-\frac{|x-y|^2}{4(t-s)}}
F(y,\2{\rho}_{i-1},\nabla\2{\rho}_{i-1},s)\,dy\,ds\quad\forall i\ge 2
\tag 1.4
$$
where
$$
v(x,t)=\left\{\aligned
&\frac{1}{(4\pi (t-T_2))^{\frac{n+2}{2}}}
\int_{\Bbb{R}^{n+2}}e^{-\frac{|x-y|^2}{4(t-T_2)}}\2{\rho}(y,T_2)\,dy
\quad\text{ for }t>T_2\\
&\2{\rho}(x,T_2)\qquad\qquad\qquad\qquad\qquad\qquad\qquad\qquad
\text{ for }t=T_2.\endaligned\right.
$$ 
Then 
$$\align
&\|v(\cdot,t)\|_{L^{\infty}(\Bbb{R}^{n+2})}
+\|\nabla v(\cdot,t)\|_{L^{\infty}(\Bbb{R}^{n+2})}\\
\le&\|\2{\rho}_1(\cdot,t)\|_{L^{\infty}(\Bbb{R}^{n+2})}
+\|\nabla\2{\rho}_1(\cdot,t)\|_{L^{\infty}(\Bbb{R}^{n+2})}\\
\le&C_1\quad\forall t\ge T_2.\tag 1.5
\endalign
$$
We claim that there exists $\delta_1\in (0,T-T_1)$ independent of $T_2$
such that 
$$
\|\2{\rho}_i(\cdot,t)\|_{L^{\infty}(\Bbb{R}^{n+2})}
+\|\nabla\2{\rho}_i(\cdot,t)\|_{L^{\infty}(\Bbb{R}^{n+2})}
\le 2C_1\quad\forall T_2\le t\le T_2+\delta_1,i\in\Bbb{Z}^+.\tag 1.6
$$
We will prove this claim by induction.
By (1.5), (1.6) holds for $i=1$. Suppose there exists $\delta_1\in (0,T-T_1)$ 
such that (1.6) holds for $i=k-1$ for some $k\ge 2$. As in \cite{LT} there 
exist a constant $C_2>0$ depending on $C_1$ and a constant
$C_3>0$ such that $\forall (x,t)\in\Bbb{R}^{n+2}\times[T_2,T_2+\delta_1)$,
$$\left\{\aligned
&|G_1(\2{\rho}_{k-1},x,t)|\le C_2\\
&|F(x,\2{\rho}_{k-1},\nabla\2{\rho}_{k-1},t)|\le C_3(2C_1)+(2C_1)^2+C_2
=C_4\, (\text{say}).\endaligned\right.\tag 1.7
$$
Then by (1.4),
$$
|\2{\rho}_k(x,t)|\le |v(x,t)|+C_4(t-T_2)\quad\forall (x,t)\in\Bbb{R}^{n+2}
\times[T_2,T_2+\delta_1)\tag 1.8
$$
and
$$\align
|\nabla\2{\rho}_k(x,t)|\le&|\nabla v(x,t)|+\int_{T_2}^t\int_{\Bbb{R}^{n+2}}
\frac{1}{(4\pi (t-s))^{\frac{n+2}{2}}}|\nabla e^{-\frac{|x-y|^2}{4(t-s)}}|
|F(y,\2{\rho}_{k-1},\nabla\2{\rho}_{k-1},s)|\,dy\,ds\\
\le&|\nabla v(x,t)|+\frac{C_4}{2}\int_{T_2}^t\int_{\Bbb{R}^{n+2}}
\frac{|x-y|}{(4\pi (t-s))^{\frac{n+2}{2}}(t-s)}
e^{-\frac{|x-y|^2}{4(t-s)}}\,dy\,ds\\
\le&|\nabla v(x,t)|+C_4C_5\int_{T_2}^t\int_{\Bbb{R}^{n+2}}
\frac{1}{(t-s)^{\frac{n+3}{2}}}e^{-\frac{|x-y|^2}{5(t-s)}}\,dy\,ds\\
\le&|\nabla v(x,t)|+C_4C_5'\sqrt{t-T_2}\quad\forall (x,t)\in\Bbb{R}^{n+2}
\times[T_2,T_2+\delta_1)\tag 1.9
\endalign
$$
for some constants $C_5>0$, $C_5'>0$, independent of $T_2$ and $\delta_1$.
Let
$$
\delta_1=\min\biggl (1,\frac{T-T_1}{2},\biggl (\frac{C_1}{C_4+C_4C_5'}
\biggr )^2\biggr ).
$$
Then by (1.5), (1.8) and (1.9), $\forall T_2
\le t\le T_2+\delta_1$,
$$
\|\2{\rho}_k(\cdot,t)\|_{L^{\infty}(\Bbb{R}^{n+2})}
+\|\nabla\2{\rho}_k(\cdot,t)\|_{L^{\infty}(\Bbb{R}^{n+2})}
\le C_1+C_4\delta_1 +C_4C_5'\sqrt{\delta_1}
\le 2C_1.
$$
Hence by induction (1.6) holds. Since
$$\align
&\2{\rho}_i(x,t)-\2{\rho}_{i-1}(x,t)\\
=&\int_{T_2}^t\int_{\Bbb{R}^{n+2}}
\frac{1}{(4\pi (t-s))^{\frac{n+2}{2}}}e^{-\frac{|x-y|^2}{4(t-s)}}
(F(y,\2{\rho}_{i-1},\nabla\2{\rho}_{i-1},s)-F(y,\2{\rho}_{i-2},
\nabla\2{\rho}_{i-2},s))\,dy\,ds\endalign
$$
for any $i\ge 3$, by (1.6) and an argument similar to the proof on P.10--11 
of \cite{LT} there exists a constant $\delta_2\in (0,\delta_1)$ independent 
of $T_2$ and depending only on $\delta_1$, $C_1$, $T_1$, and $T$ such that 
$\{\2{\rho}_i\}_{i=1}^{\infty}$ is a Cauchy sequence in $C^1(\Bbb{R}^{n+2}
\times (T_2,T_2+\delta_2))$ with norm given by 
$$
\|\psi\|=\|\psi\|_{L^{\infty}(\Bbb{R}^{n+2}\times (T_2,T_2+\delta_2))}
+\|\nabla\psi\|_{L^{\infty}(\Bbb{R}^{n+2}\times (T_2,T_2+\delta_2))}
\quad\forall\psi\in C^1(\Bbb{R}^{n+2}\times (T_2,T_2+\delta_2)).
$$ 
Let $T_2=T_1-(\delta_2/2)$. Then $T_2+\delta_2>T_1$ and there exists 
$\2{\rho}_{\infty}\in C^1(\Bbb{R}^{n+2}\times (T_2,T_2+\delta_2))$ such that
$\2{\rho}_i$ converges uniformly to $\2{\rho}_{\infty}$ in $C^1(\Bbb{R}^{n+2}
\times (T_2,T_2+\delta_2))$ as $i\to\infty$. We now extend $\2{\rho}$ beyond
the time $T_2$ by setting $\2{\rho}(x,t)=\2{\rho}_{\infty}(x,t)$ for any 
$x\in\Bbb{R}^{n+2}, T_2<t\le T_2+\delta_2$. Letting $i\to\infty$ in (1.4),
$$
\2{\rho}(x,t)=v(x,t)+\int_{T_2}^t\int_{\Bbb{R}^{n+2}}
\frac{1}{(4\pi (t-s))^{\frac{n+2}{2}}}e^{-\frac{|x-y|^2}{4(t-s)}}
F(y,\2{\rho},\nabla\2{\rho},s)\,dy\,ds\tag 1.10
$$
for any $x\in\Bbb{R}^{n+2}, T_2<t\le T_2+\delta_2$.
By Theorem 1.1 and the semi-group property of the heat equation,
$$\align
v(x,t)=&\frac{1}{(4\pi (t-T_2))^{\frac{n+2}{2}}}\int_{\Bbb{R}^{n+2}}
e^{-\frac{|x-z|^2}{4(t-T_2)}}\int_{t_0}^{T_2}\int_{\Bbb{R}^{n+2}}
\frac{1}{(4\pi (T_2-s))^{\frac{n+2}{2}}}e^{-\frac{|z-y|^2}{4(T_2-s)}}\\
&\qquad\cdot F(y,\2{\rho},\nabla\2{\rho},s)\,dy\,ds\,dz\\
=&\int_{t_0}^{T_2}\int_{\Bbb{R}^{n+2}}
\frac{1}{(4\pi (t-s))^{\frac{n+2}{2}}}e^{-\frac{|x-y|^2}{4(t-s)}}
F(y,\2{\rho},\nabla\2{\rho},s)\,dy\,ds.\tag 1.11
\endalign
$$
By (1.10) and (1.11),
$\2{\rho}$ satisfies (1.2) in $\Bbb{R}^{n+2}\times (t_0, T_2+\delta_2)$.
By (1.2) and the same argument as \cite{LT} $\2{\rho}$ is a classical 
solution of (0.11). Since $T_2+\delta_2>T_1$ and (1.6) holds, there is
a contradiction to the defintion of $T_1$. Hence $T_1=T_0$ and the theorem
follows.
\enddemo

$$
\text{Section 2}
$$

In this section we will prove various estimates for the Green function of the 
heat equation in cylindrical and punctured cylindrical domains. 
We first choose a monotone increasing function $\eta\in C^{\infty}(\Bbb{R})$, 
$0\le\eta\le 1$, such that $\eta (\tau)=0$ for any $\tau\le 1/2$ and 
$\eta (\tau)=1$ for any $\tau\ge 1$.  For any $0<\delta\le 1$, $s\in\Bbb{R}$, 
let $\eta_{\delta}(\tau)=\eta (\tau/\delta)$ and let $g_{R,\3,\delta}$ be the 
solution of 
$$\left\{\aligned
&\1_tg_{R,\3,\delta}=\Delta_xg_{R,\3,\delta}\qquad\qquad\qquad\qquad\qquad
\forall x\in B_R\setminus B_{\3}, t>s\\
&g_{R,\3,\delta}(x,t,y,s)=\Gamma (x,t,y,s)\eta ((t-s)/\delta)
\quad\forall (x,t)\in (\1 B_{\3}\cup\1 B_R)\times (s,\infty)\\
&g_{R,\3,\delta}(x,s,y,s)=0\qquad\qquad\qquad\qquad\qquad\forall x\in B_R
\endaligned\right.
$$
where $B_R\subset\Bbb{R}^m$ for some $m\in\Bbb{Z}^+$. Then by the maximum 
principle, 
$$
g_{R,\3,\delta}(x,t,y,s)\ge g_{R,\3,\delta'}(x,t,y,s)\ge 0\quad
\forall\3\le |x|,|y|\le R, 0<\delta<\delta'\le 1,t>s.\tag 2.1
$$

\proclaim{\bf Lemma 2.1}
Let $s\in\Bbb{R}$ and $\3\le |y|\le R$. Then there exists a sequence 
$\{\delta_i\}_{i=1}^{\infty}$, $\delta_i\to 0$ 
as $i\to\infty$, such that $g_{R,\3,\delta_i}(\cdot,\cdot,y,s)$ converges 
uniformly in $C^{2,1}((\2{B}_R\setminus B_{\3})\times [t_1,t_2])$ to the 
solution $g_{R,\3}(\cdot,\cdot,y,s)$ of the problem
$$\left\{\aligned
&\1_tg_{R,\3}=\Delta_xg_{R,\3}\qquad\qquad\quad\forall x\in B_R\setminus
B_{\3}, t>s\\
&g_{R,\3}(x,t,y,s)=\Gamma (x,t,y,s)\quad\forall x\in\1 B_{\3}\cup\1 B_R,t>s\\
&g_{R,\3}(x,s,y,s)=0\qquad\qquad\quad\forall x\in B_R
\endaligned\right.\tag 2.2
$$
as $i\to\infty$ for any $t_2>t_1>s$. Moreover
$$
g_{R,\3,\delta}(x,t,y,s)\le\Gamma (x,t,y,s)\quad\forall\3\le |x|,|y|\le R, 
0<\delta\le 1,t>s\tag 2.3
$$
and
$$
0\le g_{R,\3}(x,t,y,s)\le\Gamma (x,t,y,s)\quad\forall\3\le |x|,|y|\le R, t>s.
\tag 2.4
$$
\endproclaim
\demo{Proof}
Note that the result (2.4) is well-known (cf. \cite{A}). For the sake of
completeness we will give a short proof of (2.4) here. 
We will use a modification of the technique of \cite{DK} to prove (2.3).
Let $h\in C_0^{\infty}(B_R\setminus B_{\3})$ be such that $0\le h\le 1$.
For any $t>s$, let $\phi(x,\tau)$ be the solution of
$$\left\{\aligned
&\1_{\tau}\phi+\Delta\phi=0\quad\text{ in }(B_R\setminus B_{\3})\times (s,t)\\
&\phi (x,\tau)=0\qquad\text{ on }(\1 B_R\cup\1 B_{\3})\times (s,t)\\
&\phi (x,t)=h(x)\quad\text{ in }B_R\setminus B_{\3}.
\endaligned\right.
$$
By the maximum principle $0\le\phi\le 1$ on $(B_R\setminus B_{\3})\times 
(s,t)$. Hence $\1\phi/\1 n\ge 0$ on $(\1 B_R\cup\1 B_{\3})\times (s,t)$
where $\1/\1 n$ is the derivative at the boundary in the direction of the 
inward normal of the domain $B_R\setminus B_{\3}$.
Then for any $t>t_1>s$, $\3\le|y|\le R$,
$$\align
&\int_{B_R\setminus B_{\3}}(g_{R,\3,\delta}(x,t,y,s)-\Gamma (x,t,y,s))
h(x)\,dx\\
&\qquad-\int_{B_R\setminus B_{\3}}(g_{R,\3,\delta}(x,t_1,y,s)
-\Gamma (x,t_1,y,s))\phi (x,t_1)\,dx\\
=&\int_{t_1}^t\frac{\1}{\1 \tau}\biggl (\int_{B_R\setminus B_{\3}}
(g_{R,\3,\delta}(x,\tau,y,s)-\Gamma (x,\tau,y,s))\phi\,dx\biggr )\,d\tau\\
=&\int_{t_1}^t\int_{B_R\setminus B_{\3}}\phi\frac{\1}{\1\tau}
(g_{R,\3,\delta}-\Gamma)\,dx\,d\tau
+\int_{t_1}^t\int_{B_R\setminus B_{\3}}(g_{R,\3,\delta}-\Gamma)
\phi_{\tau}\,dx\,d\tau\\
=&\int_{t_1}^t\int_{B_R\setminus B_{\3}}\phi\Delta (g_{R,\3,\delta}-\Gamma)
\,dx\,d\tau
+\int_{t_1}^t\int_{B_R\setminus B_{\3}}(g_{R,\3,\delta}-\Gamma)
\phi_{\tau}\,dx\,d\tau\\
=&\int_{t_1}^t\int_{B_R\setminus B_{\3}}(g_{R,\3,\delta}-\Gamma)
(\phi_{\tau}+\Delta\phi)\,dx\,d\tau
+\int_{t_1}^t\int_{\1(B_R\setminus B_{\3})}(g_{R,\3,\delta}-\Gamma)
\frac{\1\phi}{\1 n}\,d\sigma\,d\tau.\\
\le&0.\endalign
$$
Hence
$$\align
&\int_{B_R\setminus B_{\3}}(g_{R,\3,\delta}(x,t,y,s)-\Gamma (x,t,y,s))
h(x)\,dx\\
\le&\int_{B_R\setminus B_{\3}}(g_{R,\3,\delta}(x,t_1,y,s)
-\Gamma (x,t_1,y,s))\phi (x,t_1)\,dx\quad\forall t>t_1>s\\
\Rightarrow\quad&\int_{B_R\setminus B_{\3}}(g_{R,\3,\delta}(x,t,y,s)
-\Gamma (x,t,y,s))h(x)\,dx\le 0\quad\forall t>s,\3\le|y|\le R
\quad\text{ as }t_1\searrow s.\tag 2.5
\endalign
$$
We now choose a sequence of functions $\{h_i\}_{i=1}^{\infty}\subset
C_0^{\infty}(B_R\setminus B_{\3})$, $0\le h_i\le 1$ for any $i\in\Bbb{Z}^+$,
such that $h_i(x)$ converges to $\chi_A(x)$ a.e. as $i\to\infty$ where
$A=\{x\in B_R\setminus B_{\3}:g_{R,\3,\delta}(x,t,y,s)>\Gamma (x,t,y,s)\}$. 
Putting $h=h_i$ in (2.5) and letting $i\to\infty$,
$$
\int_{B_R\setminus B_{\3}}(g_{R,\3,\delta}(x,t,y,s)
-\Gamma (x,t,y,s))_+\,dx\le 0\quad\forall t>s,\3\le|y|\le R
$$
and (2.3) follows. By (2.1), (2.3), and the Schauder estimates \cite{LSU}
for any $t_2>t_1>s$, $\3\le |y|\le R$, the sequence $\{g_{R,\3,\delta}
(\cdot,\cdot,y,s)\}_{0<\delta\le 1}$ is equi-Holder continuous on
$C^{2,1}((\2{B}_R\setminus B_{\3})\times [t_1,t_2])$. Hence by the Ascoli 
theorem and a diagonalization argument there exists a sequence 
$\{\delta_i\}_{i=1}^{\infty}$, $\delta_i\to 0$ 
as $i\to\infty$, such that $g_{R,\3,\delta_i}(\cdot,\cdot,y,s)$ converges 
uniformly in $C^{2,1}((\2{B}_R\setminus B_{\3})\times [t_1,t_2])$ to the 
solution $g_{R,\3}(\cdot,\cdot,y,s)$ of (2.2). Putting $\delta=\delta_i$
in (2.1), (2.3), and letting $i\to\infty$ we get (2.4) and the lemma follows.
\enddemo

For any $0<\delta\le 1$, $s\in\Bbb{R}$, let $g_R^{\delta}$ be the solution of 
$$\left\{\aligned
&\1_tg_R^{\delta}=\Delta_xg_R^{\delta}\qquad\qquad\qquad\qquad\qquad\quad\,\,
\forall x\in B_R, t>s\\
&g_R^{\delta}(x,t,y,s)=\Gamma (x,t,y,s)\eta ((t-s)/\delta)\quad\forall 
x\in\1 B_R,t>s\\
&g_R^{\delta}(x,s,y,s)=0\qquad\qquad\qquad\qquad\qquad\,\forall x\in B_R.
\endaligned\right.
$$
By an argument similar to the proof of (2.1) and Lemma 2.1 we have

\proclaim{\bf Lemma 2.2}
Let $s\in\Bbb{R}$ , $0<|y|\le R$, and let $\{\delta_i\}_{i=1}^{\infty}$ be 
the sequence given by Lemma 2.1. Then there exists a subsquence 
$\{\delta_i'\}_{i=1}^{\infty}$ of $\{\delta_i\}_{i=1}^{\infty}$ such that 
$g_R^{\delta_i'}(\cdot,\cdot,y,s)$ converges uniformly in $C^{2,1}
(\2{B}_R\times [t_1,t_2])$ to the solution $g_R(\cdot,\cdot,y,s)$ of (0.12) 
as $i\to\infty$ for any $t_2>t_1>s$. Moreover
$$
0\le g_R^{\delta'}(x,t,y,s)\le g_R^{\delta}(x,t,y,s)\le\Gamma (x,t,y,s)
\quad\forall 0<|x|,|y|\le R, 0<\delta\le\delta'\le 1,t>s.\tag 2.6
$$  
\endproclaim

\proclaim{\bf Lemma 2.3}
Let $m\ge 3$, $s\in\Bbb{R}$ and $0<|y|\le R$. Then there exists a sequence 
$\{\3_i\}_{i=1}^{\infty}$, $0<\3_i<R/3\quad\forall i\in\Bbb{Z}^+$ and 
$\3_i\to 0$ as $i\to\infty$, such that the sequence
$\{g_{R,\3_i}(\cdot,\cdot,y,s)\}$ converges uniformly to
$g_R(\cdot,\cdot,y,s)$ in $C^{2,1}(K)$ for any compact set $K\subset 
(\2{B}_R\setminus\{0\})\times (s,\infty)$ as $i\to\infty$.
\endproclaim
\demo{Proof}
Let $\eta$ be as before and let $\4{\eta}_{\delta}(x)=\eta (|x|/\delta)$
for any $0<\delta\le 2R/3$, $x\in\Bbb{R}^m$. Then $|\nabla\4{\eta}_{\delta}
(x)|\le C/\delta$ 
and $|\Delta\4{\eta}_{\delta}(x)|\le C/\delta^2$ on $B_R$ for some constant 
$C>0$ independent of $\delta$ and $\4{\eta}_{\delta}(x)\equiv 0$ for any 
$|x|\le\delta/2$. Let $h\in C_0^{\infty}(B_R)$ be such that $0\le h\le 1$. 
For any $t>s$, let $\phi(x,\tau)$ be the solution of
$$\left\{\aligned
&\1_{\tau}\phi+\Delta\phi=0\quad\text{ in }B_R\times (s,t)\\
&\phi (x,\tau)=0\qquad\text{ on }\1 B_R\times (s,t)\\
&\phi (x,t)=h(x)\quad\text{ in }B_R.
\endaligned\right.\tag 2.7
$$
Now by (2.1), (2.3), and (2.6), for any $0<\3\le R/3$, $0<\delta\le 1$, $t>s$,
$$\align
&\int_{B_R\setminus B_{\3}}(g_R^{\delta}(x,t,y,s)-g_{R,\3,\delta}(x,t,y,s)
h(x)\4{\eta}_{2\3}(x)\,dx\\
=&\int_s^t\frac{\1}{\1 \tau}\biggl (\int_{B_R\setminus B_{\3}}
(g_R^{\delta}(x,\tau,y,s)-g_{R,\3,\delta}(x,\tau,y,s))\phi\4{\eta}_{2\3}
\,dx\biggr )\,d\tau\\
=&\int_s^t\int_{B_R\setminus B_{\3}}\phi\4{\eta}_{2\3}\frac{\1}{\1\tau}
(g_R^{\delta}-g_{R,\3,\delta})\,dx\,d\tau
+\int_s^t\int_{B_R\setminus B_{\3}}(g_R^{\delta}-g_{R,\3,\delta})
\phi_{\tau}\4{\eta}_{2\3}\,dx\,d\tau\\
=&\int_s^t\int_{B_R\setminus B_{\3}}\phi\4{\eta}_{2\3}\Delta 
(g_R^{\delta}-g_{R,\3,\delta})\,dx\,d\tau
+\int_s^t\int_{B_R\setminus B_{\3}}(g_R^{\delta}-g_{R,\3,\delta})
\phi_{\tau}\4{\eta}_{2\3}\,dx\,d\tau\\
=&\int_s^t\int_{B_R\setminus B_{\3}}(g_R^{\delta}-g_{R,\3,\delta})
\Delta (\phi\4{\eta}_{2\3})\,dx\,d\tau
+\int_s^t\int_{B_R\setminus B_{\3}}(g_R^{\delta}-g_{R,\3,\delta})
\phi_{\tau}\4{\eta}_{2\3}\,dx\,d\tau\\
=&\int_s^t\int_{B_R\setminus B_{\3}}(g_R^{\delta}-g_{R,\3,\delta})
(\phi_{\tau}+\Delta\phi)\4{\eta}_{2\3}\,dx\,d\tau\\
&\qquad+\int_s^t\int_{\3\le |x|\le 2\3}(g_R^{\delta}-g_{R,\3,\delta})
(2\nabla\phi\cdot\nabla\4{\eta}_{2\3}+\phi\Delta\4{\eta}_{2\3})
\,dx\,d\tau\\
\le &\frac{C}{\3^2}\int_s^t\int_{\3\le |x|\le 2\3}\Gamma (x,\tau,y,s)
\,dx\,d\tau
+\frac{C}{\3}\int_s^t\int_{\3\le |x|\le 2\3}\Gamma (x,\tau,y,s)|\nabla\phi|
\,dx\,d\tau\\
=&I_1+I_2\tag 2.8
\endalign
$$
for some constant $C>0$ independent of $0<\3\le R/3$ and $0<\delta\le 1$.
Let $0<\3\le\min (|y|/10,R/3)$. Then for any $\3\le |x|\le 2\3$, 
$$
|x|\le\frac{|y|}{5}\quad\Rightarrow\quad |x-y|\ge |y|-|x|\ge\frac{4}{5}|y|.
$$
Hence 
$$\align
\Gamma (x,\tau,y,s)\le&\frac{1}{(4\pi (\tau -s))^{\frac{m}{2}}}
e^{-\frac{4|y|^2}{25(\tau -s)}}
=\biggl (\frac{|y|}{(4\pi (\tau -s))^{\frac{1}{2}}}
\biggr )^me^{-\frac{2|y|^2}{25(\tau -s)}}\cdot
\frac{e^{-\frac{2|y|^2}{25(\tau -s)}}}{|y|^m}\\
\le&\frac{C}{|y|^m}e^{-\frac{2|y|^2}{25(\tau -s)}}
\quad\forall \3\le |x|\le 2\3, \tau>s.\tag 2.9
\endalign
$$
Thus
$$
I_1\le C\frac{\3^{m-2}}{|y|^m}e^{-\frac{2|y|^2}{25(t-s)}}(t-s).\tag 2.10
$$
Now by (2.7),
$$\align
&\frac{\1}{\1\tau}\int_{B_R}|\nabla\phi|^2\,dx=2\int_{B_R}\nabla\phi\cdot
\nabla\phi_{\tau}\,dx=-2\int_{B_R}\Delta\phi\cdot
\phi_{\tau}\,dx=2\int_{B_R}(\Delta\phi)^2\,dx\\
\Rightarrow&\quad \sup_{s\le\tau\le t}\int_{B_R}|\nabla\phi|^2\,dx
\le\int_{B_R}|\nabla h|^2\,dx.\tag 2.11
\endalign
$$
By (2.9) and (2.11),
$$\align
I_2\le&\frac{C}{\3}\cdot\frac{e^{-\frac{2|y|^2}{25(t-s)}}}{|y|^m}
\int_s^t\int_{B_R}|\nabla\phi|\,dx\,d\tau\\
\le&C\3^{\frac{m}{2}-1}\cdot\frac{e^{-\frac{2|y|^2}{25(t-s)}}}{|y|^m}
\int_s^t\biggl (\int_{B_R}|\nabla\phi|^2\,dx\biggr )^{\frac{1}{2}}\,d\tau\\
\le&C\3^{\frac{m}{2}-1}\cdot\frac{e^{-\frac{2|y|^2}{25(t-s)}}}{|y|^m}
\biggl (\int_{B_R}|\nabla h|^2\,dx\biggr )^{\frac{1}{2}}(t-s).\tag 2.12
\endalign
$$
By (2.8), (2.10), and (2.12),
$$\align
&\int_{B_R\setminus B_{\3}}(g_R^{\delta}(x,t,y,s)-g_{R,\3,\delta}(x,t,y,s)
h(x)\4{\eta}_{2\3}(x)\,dx\\
\le&C\frac{\3^{m-2}}{|y|^m}e^{-\frac{2|y|^2}{25(t-s)}}(t-s)
+C\3^{\frac{m}{2}-1}\cdot\frac{e^{-\frac{2|y|^2}{25(t-s)}}}{|y|^m}
\biggl (\int_{B_R}|\nabla h|^2\,dx\biggr )^{\frac{1}{2}}(t-s)\quad\forall
t>s.\tag 2.13
\endalign
$$
Let $\{\delta_i'\}_{i=1}^{\infty}$ be the sequence given by Lemma 2.2.
Putting $\delta=\delta_i'$ in (2.13) and letting $i\to\infty$ 
by Lemma 2.1, Lemma 2.2, and the Lebesgue dominated convergence theorem,
$$\align
&\int_{B_R\setminus B_{\3}}
(g_R(x,t,y,s)-g_{R,\3}(x,t,y,s)h(x)\4{\eta}_{2\3}(x)\,dx\\
\le&C\frac{\3^{m-2}}{|y|^m}e^{-\frac{2|y|^2}{25(t-s)}}(t-s)
+C\3^{\frac{m}{2}-1}\cdot\frac{e^{-\frac{2|y|^2}{25(t-s)}}}{|y|^m}
\biggl (\int_{B_R}|\nabla h|^2\,dx\biggr )^{\frac{1}{2}}(t-s)
\quad\forall t>s.\tag 2.14
\endalign
$$
By (2.4) and the Schauder estimates \cite{LSU} the sequence 
$\{g_{R,\3}(x,t,y,s)\}_{0<\3\le R/3}$ is equi-Holder continuous
in $C^{2,1}(K)$ for any compact subset $K\subset (\2{B}_R\setminus\{0\})
\times [t_1,t_2]$ for any $t_2>t_1>s$. By the Ascoli theorem and a 
diagonalization argument 
there exists a sequence $\{\3_i\}_{i=1}^{\infty}$, $0<\3_i<R/3\quad\forall 
i\in\Bbb{Z}^+$ and $\3_i\to 0$ as $i\to\infty$, such that the sequence
$\{g_{R,\3_i}(\cdot,\cdot,y,s)\}$ converges uniformly to some function
$\4{g}_R(\cdot,\cdot,y,s)$ on every compact subset of $(\2{B}_R\setminus\{0\})
\times (s,\infty)$ as $i\to\infty$.

Putting $\3=\3_i$ in (2.14) and letting $i\to\infty$,  
$$
\int_{B_R\setminus\{0\}}(g_R(x,t,y,s)-\4{g}_R(x,t,y,s))h(x)\,dx\le 0
\quad\forall t>s, 0<|y|\le R.\tag 2.15
$$
We now choose a sequence of functions $\{h_i\}_{i=1}^{\infty}\subset
C_0^{\infty}(B_R)$, $0\le h_i\le 1$ for any $i\in\Bbb{Z}^+$,
such that $h_i$ converges to $\chi_A$ a.e. as $i\to\infty$ where
$A=\{x\in B_R:g_R(x,t,y,s)>\4{g}_R(x,t,y,s)\}$. 
Putting $h=h_i$ in (2.15) and letting $i\to\infty$,
$$\align
&\int_{B_R\setminus\{0\}}(g_R(x,t,y,s)-\4{g}_R(x,t,y,s))_+\,dx\le 0
\quad\forall t>s, 0<|y|\le R\\
\Rightarrow\quad&g_R(x,t,y,s)\le\4{g}_R(x,t,y,s)\qquad\qquad\qquad
\qquad\quad\forall 0<|x|,|y|\le R, t>s.\tag 2.16
\endalign
$$
Interchanging the role of $g_R(x,t,y,s)$ and $\4{g}_R(x,t,y,s)$ and repeating
the above argument we get
$$
g_R(x,t,y,s)\ge\4{g}_R(x,t,y,s)\quad\forall 0<|x|,|y|\le R, t>s.\tag 2.17
$$
By (2.16) and (2.17),
$$
\4{g}_R(x,t,y,s)=g_R(x,t,y,s)\quad\forall 0<|x|,|y|\le R, t>s\tag 2.18
$$
and the lemma follows.
\enddemo

For any $0<\3<R$, let $G_{R,\3}=G_{R,\3}(x,t,y,s)$ be the Green function
for the heat equation in $(B_R\setminus B_{\3})\times (-\infty,\infty)$.
Then 
$$
G_{R,\3}(x,t,y,s)=\Gamma(x,t,y,s)-g_{R,\3}(x,t,y,s)\quad\forall \3\le |x|,
|y|\le R,t>s.\tag 2.19
$$
Then by Lemma 2.3 and the uniqueness of the Green function for the heat 
equation we have the following corollary.

\proclaim{\bf Corollary 2.4}
Let $m\ge 3$, $s\in\Bbb{R}$ and $0<|y|\le R$. Then $G_{R,\3}(\cdot,\cdot,y,s)$ 
will converge uniformly to $G_R(\cdot,\cdot,y,s)$ in $C^{2,1}(K)$ for any 
compact set $K\subset (\2{B}_R\setminus\{0\})\times (s,\infty)$ as $\3\to 0$.
\endproclaim

We now let $G_{R,\3}^{\ast}=G_{R,\3}^{\ast}(y,s,x,t)$ and
$G_R^{\ast}=G_R^{\ast}(y,s,x,t)$, $s<t$,  be the Green function for the 
adjoint heat equation $\1_su+\Delta u=0$ in $(B_R\setminus B_{\3})\times 
(-\infty,\infty)$ and $B_R\times (-\infty,\infty)$ respectively.
Then by a similar argument as the proof of Lemma 2.3 and Corollary 2.4
we have the following result.

\proclaim{\bf Corollary 2.5}
Let $m\ge 3$. Then for any $t\in\Bbb{R}$, $0<|x|\le R$, 
$G_{R,\3}^{\ast}(\cdot,\cdot,x,t)$ will converge uniformly  to 
$G_R^{\ast}(\cdot,\cdot,x,t)$ in $C^{2,1}(K)$ for any compact set 
$K\subset(\2{B}_R\setminus\{0\})\times (-\infty,t)$ as $\3\to 0$.  
\endproclaim

Now by \cite{F}, 
$$
G_{R,\3}(x,t,y,s)=G_{R,\3}^{\ast}(y,s,x,t)\quad\forall \3\le |x|,|y|\le R, 
t>s\tag 2.20
$$
and
$$
G_R(x,t,y,s)=G_R^{\ast}(y,s,x,t)\quad\forall |x|,|y|\le R, t>s.\tag 2.21
$$
Hence by (2.20), (2.21), and Corollary 2.5 we have the following result.

\proclaim{\bf Corollary 2.6}
Let $m\ge 3$. Then for any $t\in\Bbb{R}$, $0<|x|\le R$, 
$G_{R,\3}(x,t,\cdot,\cdot)$ will converge uniformly to $G_R(x,t,\cdot,\cdot)$ 
in $C^{2,1}(K)$ for any compact set $K\subset(\2{B}_R\setminus\{0\})\times 
(-\infty,t)$ as $\3\to 0$.
\endproclaim

Let $G_{\infty,1}(x,t,y,s)$ be the Green function for the heat equation in the
domain $(\Bbb{R}^{m-1}\setminus\2{B}_1)\times (-\infty,\infty)$.
By an argument similar to the proof of (ii) of Lemma 1.3 of \cite{Hu} we
have

\proclaim{\bf Lemma 2.7}
For any $T>0$, there exist constants $C>0$ and $c>0$ such that
$$
0\le G_{\infty,1}(x,t,y,s)\le C\frac{\delta_1(y)}{(t-s)^{\frac{m+1}{2}}}
e^{-c|x-y|^2/(t-s)}\quad\forall |x|,|y|>1,0\le s<t\le T
$$
and
$$
0\le G_R(x,t,y,s)\le C\frac{\delta_R(y)}{(t-s)^{\frac{m+1}{2}}}
e^{-c|x-y|^2/(t-s)}\quad\forall |x|,|y|<R,0\le s<t\le T, R>0
$$
holds where $\delta_1(y)=$dist$(y,\1 B_1)$ and $\delta_R(y)
=$dist$(y,\1 B_R)$ respectively.
\endproclaim

\proclaim{\bf Corollary 2.8}
For any $T>0$ there exist constants $C>0$ and $c>0$ such that
$$
0\le\frac{\1 G_{\infty,1}}{\1 n_y}(x,t,y,s)
\le\frac{C}{(t-s)^{\frac{m+1}{2}}}e^{-c|x-y|^2/(t-s)}\quad\forall 
|x|>1,|y|=1,0\le s<t\le T
$$
holds where $\1/\1 n_y$ is the derivative in the direction of the inward normal
$n_y=y/|y|$ at the point $y\in\1 B_1(0)$ with respect to the domain 
$\Bbb{R}^m\setminus\2{B}_1$.
\endproclaim
 
\proclaim{\bf Corollary 2.9}
For any $T>0$ there exist constants $C>0$ and $c>0$ such that
$$
0\le\frac{\1 G_R}{\1 n_y}(x,t,y,s)
\le\frac{C}{(t-s)^{\frac{m+1}{2}}}
e^{-c|x-y|^2/(t-s)}\quad\forall |x|<R,|y|=R,0\le s<t\le T, R>0
$$
holds where $\1/\1 n_y$ is the derivative in the direction of the inward
normal $n_y=-y/|y|$ at the point $y\in\1 B_R$ with respect to the domain 
$B_R$.
\endproclaim

\proclaim{\bf Corollary 2.10}
For any $T>0$ there exist constants $C>0$ and $c>0$ such that
$$
0\le\frac{\1 G_{R,\3}}{\1 n_y}(x,t,y,s)
\le\frac{C}{(t-s)^{\frac{m+1}{2}}}e^{-c|x-y|^2/(t-s)}\quad\forall |y|=\3,
\3<|x|<R,0\le s<t\le T\tag 2.22
$$
holds where $\1/\1 n_y$ is the derivative in the direction of the inward
normal $n_y=y/|y|$ at the point $y\in\1 B_{\3}$ with respect to the 
domain $B_R\setminus\2{B}_{\3}$.
\endproclaim
\demo{Proof}
By scaling,
$$
G_{R,\3}(x,t,y,s)=\3^{-m}G_{R/\3,1}(x/\3,t/\3^2,y/\3,s/\3^2).
$$
Hence
$$
\frac{\1 G_{R,\3}}{\1 n_y}(x,t,y,s)=\3^{-m-1}\frac{\1 G_{R/\3,1}}{\1 n_{y'}}
(x/\3,t/\3^2,y/\3,s/\3^2),y'=y/\3.\tag 2.23
$$
By the maximum principle (cf. \cite{F}), $\forall R_2>R_1>1$,
$$\align
&0\le G_{R_1,1}(x,t,y,s)\le G_{R_2,1}(x,t,y,s)\le G_{\infty,1}(x,t,y,s)
\quad\forall 1<|x|,|y|\le R_1,t>s\\
\Rightarrow\quad&0\le\frac{\1 G_{R_1,1}}{\1 n_y}(x,t,y,s)
\le\frac{\1 G_{\infty,1}}{\1 n_y}(x,t,y,s)\quad\forall 1<|x|<R_1,|y|=1,t>s, 
R_1>1.\tag 2.24
\endalign
$$
By (2.23), (2.24), and Corollary 2.8 we get (2.22) and the lemma follows.
\enddemo

By an argument similar to the proof of Corollary 2.10 but with $G_R$ and
Corollary 2.9 replacing $G_{\infty,1}$ and Corollary 2.8 in the proof
we have

\proclaim{\bf Corollary 2.11}
For any $T>0$ there exist constants $C>0$ and $c>0$ such that
$$
0\le\frac{\1 G_{R,\3}}{\1 n_y}(x,t,y,s)
\le\frac{C}{(t-s)^{\frac{m+1}{2}}}e^{-c|x-y|^2/(t-s)}\quad\forall |y|=R,
\3<|x|<R,0\le s<t\le T\tag 2.25
$$
holds where $\1/\1 n_y$ is the derivative in the direction of the inward
normal $n_y=-y/|y|$ at the point $y\in\1 B_R$ with respect to the 
domain $B_R\setminus\2{B}_{\3}$.
\endproclaim

$$
\text{Section 3}
$$

In this section we will use the Green function estimates obtained in section
two to prove that under an uniform boundedness condition on a solution of 
(0.11) in $(\Bbb{R}^{n+2}\setminus\{0\})\times (t_0,T_0)$ the solution has 
removable singularities on the line $\{0\}\times (t_0,T_0)$. We will also 
prove the uniqueness of solution of (0.8). 

\proclaim{\bf Lemma 3.1}
Let $T'\in (t_0,T]$ and let $\2{\rho}$ be a solution of (0.11)  in 
$(\Bbb{R}^{n+2}\setminus\{0\})\times (t_0,T')$ which satisfies (1.1) and 
$$
\sup_{t_0\le t<T'}(\|\2{\rho}(\cdot,t)\|_{L^{\infty}(\Bbb{R}^{n+2}\setminus
\{0\})}
+\|\nabla\2{\rho}(\cdot,t)\|_{L^{\infty}(\Bbb{R}^{n+2}\setminus
\{0\})})<\infty.\tag 3.1
$$
Then $\2{\rho}$ can be extended to a solution of (0.11) in 
$\Bbb{R}^{n+2}\times (t_0,T')$ and $\2{\rho}$ satisfies (1.2) and (1.3) in 
$\Bbb{R}^{n+2}\times (t_0,T')$.
\endproclaim
\demo{Proof}
Let $|x|>0$. By standard parabolic theory (cf. \cite{A} and \cite{LSU}), for 
any $t_0<t<T'$ and $R>\3>0$ such that $\3\le |x|\le R$,
$$\align
&\2{\rho}(x,t)\\
=&\int_{t_0}^t\int_{B_R\setminus B_{\3}}G_{R,\3}(x,t,y,s)
F(y,\2{\rho},\nabla\2{\rho},s)\,dy\,ds
+\int_{t_0}^t\int_{\1 B_{\3}}\frac{\1 G_{R,\3}}{\1 n_y}(x,t,y,s)
\2{\rho}(y,s)\,d\sigma (y)\,ds\\
&\qquad+\int_{t_0}^t\int_{\1 B_R}\frac{\1 G_{R,\3}}{\1 n_y}(x,t,y,s)
\2{\rho}(y,s)\,d\sigma (y)\,ds\\
=&I_{1,R}^{\3}+I_{2,R}^{\3}+I_{3,R}^{\3}\tag 3.2
\endalign
$$
where $\1/\1 n_y$ is the derivative in the direction of the inward normal
of the domain $B_R\setminus B_{\3}$ at $y\in\1 B_R\cup\1 B_{\3}$. Since 
$\2{\rho}$ satisfies (3.1), by an argument similar to the proof on P.10
of \cite{LT} there exists a constant $C>0$ such that
$$
|F(x,\2{\rho},\nabla\2{\rho},t)|\le C\quad\forall x\in\Bbb{R}^{n+2}\setminus
\{0\},t_0\le t<T'.\tag 3.3
$$
By (2.4), (2.19), (3.3), Corollary 2.6, and Lebesgue dominated 
convergence theorem,
$$
\lim_{\3\to 0}I_{1,R}^{\3}=\int_{t_0}^t\int_{B_R\setminus\{0\}}G_R(x,t,y,s)
F(y,\2{\rho},\nabla\2{\rho},s)\,dy\,ds.\tag 3.4
$$
By Corollary 2.10 (2.22) holds with $m=n+2$. Hence
$$\align
&|I_{2,R}^{\3}|\le C\int_{t_0}^t\int_{\1 B_{\3}}\frac{1}{(t-s)^{\frac{n+3}{2}}}
e^{-c|x-y|^2/(t-s)}\,d\sigma (y)\,ds\le\frac{C\3^{n+1}}{(|x|-\3)^{n+3}}\\
\Rightarrow\quad&\lim_{\3\to 0}I_{2,R}^{\3}=0.\tag 3.5
\endalign
$$
Now
$$\align
&I_{3,R}^{\3}\\
=&\int_{t_0}^{t-\delta}\int_{\1 B_R}
\frac{\1 G_{R,\3}}{\1 n_y}(x,t,y,s)\2{\rho}(y,s)\,d\sigma (y)\,ds
+\int_{t-\delta}^t\int_{\1 B_R}\frac{\1 G_{R,\3}}{\1 n_y}(x,t,y,s)
\2{\rho}(y,s)\,d\sigma (y)\,ds\\
=&J_{1,\delta}^{\3}+J_{2,\delta}^{\3}\quad\forall 0<\delta<t-t_0.\tag 3.6
\endalign
$$
By (3.1) and Corollary 2.6,
$$
\lim_{\3\to 0}J_{1,\delta}^{\3}=\int_{t_0}^{t-\delta}\int_{\1 B_R}
\frac{\1 G_R}{\1 n_y}(x,t,y,s)\2{\rho}(y,s)\,d\sigma (y)\,ds
\quad\forall 0<\delta<t-t_0.\tag 3.7
$$
By Corollary 2.11 (2.25) holds with $m=n+2$. Then
$$\align
|J_{2,\delta}^{\3}|\le&C\int_{t-\delta}^t\int_{\1 B_R}
\frac{1}{(t-s)^{\frac{n+3}{2}}}e^{-c|x-y|^2/(t-s)}
\,d\sigma (y)\,ds\le\frac{C\delta}{(R-|x|)^{n+3}}.\tag 3.8
\endalign
$$
Similarly by Corollary 2.9,
$$
\biggl|\int_{t-\delta}^t\int_{\1 B_R}\frac{\1 G_R}{\1 n_y}(x,t,y,s)
\2{\rho}(y,s)\,d\sigma (y)\,ds\biggr|
\le\frac{C\delta}{(R-|x|)^{n+3}}.\tag 3.9
$$
Now
$$\align
&\biggl|I_{3,R}^{\3}-\int_{t_0}^t\int_{\1 B_R}\frac{\1 G_R}{\1 n_y}(x,t,y,s)
\2{\rho}(y,s)\,d\sigma (y)\,ds\biggr|\\
\le&\biggl|J_{1,\delta}^{\3}-\int_{t_0}^{t-\delta}\int_{\1 B_R}
\frac{\1 G_R}{\1 n_y}(x,t,y,s)\2{\rho}(y,s)\,d\sigma (y)\,ds\biggr|
+|J_{2,\delta}^{\3}|\\
&\qquad +\biggl|\int_{t-\delta}^t\int_{\1 B_R}\frac{\1 G_R}{\1 n_y}
(x,t,y,s)\2{\rho}(y,s)\,d\sigma (y)\,ds\biggr|.\tag 3.10
\endalign
$$
Letting $\3\to 0$ in (3.10), by (3.7), (3.8), and (3.9),
$$\align
&\limsup_{\3\to 0}\biggl|I_{3,R}^{\3}-\int_{t_0}^t\int_{\1 B_R}
\frac{\1 G_R}{\1 n_y}(x,t,y,s)\2{\rho}(y,s)\,d\sigma (y)\,ds\biggr|
\le\frac{C\delta}{(R-|x|)^{n+3}}\quad\forall 0<\delta<t-t_0\\
\Rightarrow\quad&\lim_{\3\to 0}\biggl|I_{3,R}^{\3}-\int_{t_0}^t\int_{\1 B_R}
\frac{\1 G_R}{\1 n_y}(x,t,y,s)\2{\rho}(y,s)\,d\sigma (y)\,ds\biggr|=0
\quad\text{ as }\delta\to 0.\tag 3.11
\endalign
$$
Thus letting $\3\to 0$ in (3.2), by (3.4), (3.5), and (3.11),
$\forall 0<|x|<R$,
$$\align
&\2{\rho}(x,t)\\
=&\int_{t_0}^t\int_{B_R\setminus\{0\}}G_R(x,t,y,s)
F(y,\2{\rho},\nabla\2{\rho},s)\,dy\,ds
+\int_{t_0}^t\int_{\1 B_R}\frac{\1 G_R}{\1 n_y}(x,t,y,s)
\2{\rho}(y,s)\,d\sigma (y)\,ds.\tag 3.12
\endalign
$$
Since $G_R(x,t,y,s)$ increases to $\Gamma (x,t;y,s)$ as $R\to\infty$, by 
(0.13), (3.3),  and the Lebesgue dominated convergence theorem,
$$\align
&\lim_{R\to\infty}\int_{t_0}^t\int_{B_R\setminus\{0\}}G_R(x,t,y,s)
F(y,\2{\rho},\nabla\2{\rho},s)\,dy\,ds\\
=&\int_{t_0}^t\int_{\Bbb{R}^{n+2}\setminus\{0\}}\Gamma (x,t,y,s)
F(y,\2{\rho},\nabla\2{\rho},s)\,dy\,ds.\tag 3.13
\endalign
$$
By Corollary 2.9,
$$\align
\biggl |\int_{t_0}^t\int_{\1 B_R}\frac{\1 G_R}{\1 n_y}(x,t,y,s)\2{\rho}(y,s)
\,d\sigma (y)\,ds\biggr |
\le&C\int_{t_0}^t\int_{\1 B_R}\frac{1}{(t-s)^{\frac{n+3}{2}}}
e^{-c|x-y|^2/(t-s)}\,d\sigma (y)\,ds\\
\le&\frac{CR^{n+1}}{(R-|x|)^{n+3}}.\tag 3.14
\endalign
$$
Letting $R\to\infty$ in (3.12), by (3.13) and (3.14) we get that $\2{\rho}$ 
satisfies 
$$
\2{\rho}(x,t)=\int_{t_0}^t\int_{\Bbb{R}^{n+2}\setminus\{0\}}
\frac{1}{(4\pi (t-s))^{\frac{n+2}{2}}}e^{-\frac{|x-y|^2}{4(t-s)}}
F(y,\2{\rho},\nabla\2{\rho},s)\,dy\,ds\tag 3.15
$$
$\forall (x,t)\in (\Bbb{R}^{n+2}\setminus\{0\})\times [t_0,T')$.
Since the right hand side of (3.15) is a continuous function on 
$\Bbb{R}^{n+2}\times (t_0,T')$, we can extend $\2{\rho}(x,t)$ to a 
function on $\Bbb{R}^{n+2}\times (t_0,T')$ by letting $\2{\rho}(x,t)$ 
equal to the right hand side of (3.15) for $x=0$, $t_0<t<T'$. 
Then the extended function $\2{\rho}(x,t)$
is a continuous function on $\Bbb{R}^{n+2}\times (t_0,T')$ and has
continuous first derivatives in $x$ on $\Bbb{R}^{n+2}\times (t_0,T')$.
Hence $\2{\rho}(x,t)$ satisfies (1.2) on $\Bbb{R}^{n+2}\times (t_0,T')$.
By the same argument as that on P.11--12 of \cite{LT} $\2{\rho}$ is a 
classical solution of (0.11) in $\Bbb{R}^{n+2}\times (t_0,T')$. By (3.1)
$\2{\rho}$ satisfies (1.3) and the lemma follows.
\enddemo

\proclaim{\bf Theorem 3.2}
Let $T'\in (t_0,T]$. Suppose $\2{\rho}_1$ and $\2{\rho}_2$ are solutions of 
(0.11) and (1.1) in $(\Bbb{R}^{n+2}\setminus\{0\})\times (t_0,T')$ which
 satisfies (3.1). Then $\2{\rho}_1\equiv\2{\rho}_2$ in $(\Bbb{R}^{n+2}
\setminus\{0\})\times (t_0,T')$.
\endproclaim
\demo{Proof}
Let $T_1\ge t_0$ be the maximal time such that $\2{\rho}_1\equiv\2{\rho}_2$ 
in $(\Bbb{R}^{n+2}\setminus\{0\})\times (t_0,T_1)$. Suppose $T_1<T'$. 
By Lemma 3.1 both $\2{\rho}_1$ and $\2{\rho}_2$ can be extended to solutions
of (0.11) in $\Bbb{R}^{n+2}\times (t_0,T')$ and both satisfy (1.2) and (1.3) 
in $\Bbb{R}^{n+2}\times (t_0,T')$. Then by (1.2)
$$\align
&\2{\rho}_1(x,t)-\2{\rho}_2(x,t)\\
=&\int_{T_1}^t\int_{\Bbb{R}^{n+2}}
\frac{1}{(4\pi (t-s))^{\frac{n+2}{2}}}e^{-\frac{|x-y|^2}{4(t-s)}}
[F(y,\2{\rho}_1,\nabla\2{\rho}_1,s)-F(y,\2{\rho}_2,\nabla\2{\rho}_2,s)]
\,dy\,ds\tag 3.16
\endalign
$$
Hence
$$
\align
&|\nabla\2{\rho}_1(x,t)-\nabla\2{\rho}_2(x,t)|\\
\le&\int_{T_1}^t\int_{\Bbb{R}^{n+2}}\frac{|x-y|}{2(t-s)}
\frac{1}{(4\pi (t-s))^{\frac{n+2}{2}}}e^{-\frac{|x-y|^2}{4(t-s)}}
|F(y,\2{\rho}_1,\nabla\2{\rho}_1,s)-F(y,\2{\rho}_2,\nabla\2{\rho}_2,s)|
\,dy\,ds\\
\le&\int_{T_1}^t\int_{\Bbb{R}^{n+2}}\frac{C}{(t-s)^{\frac{n+3}{2}}}
e^{-\frac{|x-y|^2}{5(t-s)}}
|F(y,\2{\rho}_1,\nabla\2{\rho}_1,s)-F(y,\2{\rho}_2,\nabla\2{\rho}_2,s)|
\,dy\,ds\tag 3.17
\endalign
$$
By (1.3) and the argument on P.10--11 of \cite{LT} there exists a constant
$C>0$ such that
$$
|F(y,\2{\rho}_1,\nabla\2{\rho}_1,s)-F(y,\2{\rho}_2,\nabla\2{\rho}_2,s)|
\le C(|\2{\rho}_1-\2{\rho}_2|+|\nabla\2{\rho}_1-\nabla\2{\rho}_2|)
\quad\forall y\in\Bbb{R}^{n+2},s\in [t_0,T'].
\tag 3.18
$$
Let
$$
E(t)=\sup_{T_1\le s<t}
(\|\2{\rho}_1(\cdot,s)-\2{\rho}_2(\cdot,s)\|_{L^{\infty}(\Bbb{R}^{n+2})}
+\|\nabla\2{\rho}_1(\cdot,s)-\nabla\2{\rho}_2(\cdot,s)\|_{L^{\infty}
(\Bbb{R}^{n+2})}).
$$
By (3.16), (3.17) and (3.18), there exist constants $C_5>0$ and $C_6>0$ such 
that
$$
E(t)\le C_5(t-T_1)E(t)+C_6\sqrt{t-T_1}E(t)\quad\forall T_1\le t<T'.
\tag 3.19
$$
Let 
$$
\delta=\min(1,1/(4(C_5+C_6)^2)).
$$
Then by (3.19) for any $T_1\le t\le T_1+\delta$,
$$\align
E(t)\le\frac{E(t)}{2}\quad&\Rightarrow\quad 
E(t)=0\quad\forall T_1\le t\le T_1+\delta\\
&\Rightarrow\quad\2{\rho}_1(x,t)=\2{\rho}_2(x,t)\quad\forall x\in
\Bbb{R}^{n+2},T_1\le t\le T_1+\delta.\endalign
$$
This contradicts the choice of $T_1$. Hence $T_1=T'$ and the theorem follows.
\enddemo

By (0.9) and Theorem 3.2 we have the following uniqueness result.

\proclaim{\bf Theorem 3.3}
Let $0<t_0<T_0<T$. Suppose $\4{\rho}_1$, $\4{\rho}_2$, are solutions of
(0.8) in $\Bbb{R}^+\times (t_0,T_0)$ satisfying
$$
\4{\rho}_i(r,t_0)=0\quad\forall r\ge 0,i=1,2
$$
and
$$
\sup_{t_0\le t<T'}(\|\widetilde{\rho}_i(\cdot ,t)\|_{L^{\infty}(\Bbb{R}^+)}
+\|\1\widetilde{\rho}_i/\1 r(\cdot ,t)\|_{L^{\infty}(\Bbb{R}^+)})
<\infty\quad\forall i=1,2
$$
for any $T'\in (t_0,T_0)$. Then
$$
\4{\rho}_1(r,t)=\4{\rho}_2(r,t)\quad\forall r\ge 0, t_0\le t<T_0.
$$
\endproclaim

$$
\text{Section 4}
$$

In this section we will prove the removable singularities property of  
the solution of the heat equation.

\proclaim{\bf Theorem 4.1}
Let $m\ge 3$ and let $\Omega\subset\Bbb{R}^m$ be a domain. Suppose $u$ is a 
solution of the heat equation in $(\Omega\setminus\{0\})\times (0,T)$. Then 
$u$ has removable singularities at $\{0\}\times (0,T)$ if and only if
there exists $\2{B}_{\delta}\subset\Omega$ such that (0.7) holds.
\endproclaim
\demo{Proof}
If $u$ has removable singularities at $\{0\}\times (0,T)$, then there exists
a solution $v$ of the heat equation in $\Omega\times (0,T)$ such that
$u=v$ on $(\Omega\setminus\{0\})\times (0,T)$. Choose $\delta>0$ such that
$\2{B}_{\delta}\subset\Omega$. Then (0.7) holds for $v$. Hence (0.7)
holds for $u$.

Now suppose that there exists $\2{B}_R\subset\Omega$ such that 
(0.7) holds for $0<t_1<t_2<T$. By standard parabolic theory,
$\forall 0<\3<|x|<R$, $t_1<t\le t_2$,
$$\align
&u(x,t)\\
=&\int_{B_R\setminus B_{\3}}G_{R,\3}(x,t,y,t_1)u(y,t_1)\,dy
+\int_{t_1}^t\int_{\1 B_{\3}}\frac{\1 G_{R,\3}}{\1 n_y}(x,t,y,s)
u(y,s)\,d\sigma (y)\,ds\\
&\qquad+\int_{t_1}^t\int_{\1 B_R}\frac{\1 G_{R,\3}}{\1 n_y}(x,t,y,s)
u(y,s)\,d\sigma (y)\,ds\\
=&I_{1,\3}+I_{2,\3}+I_{3,\3}\tag 4.1
\endalign
$$
where $\1/\1 n_y$ is the derivative in the direction of the inward normal
of the domain $B_R\setminus B_{\3}$ at $y\in\1 B_R\cup\1 B_{\3}$.
By an argument similar to the proof of Lemma 3.1 we get
$$
\lim_{\3\to 0}I_{1,\3}=\int_{B_R\setminus\{0\}}G_R(x,t,y,t_1)u(y,t_1)
\,dy\quad\forall 0<|x|<R, t_1<t\le t_2\tag 4.2
$$
and
$$
\lim_{\3\to 0}I_{3,\3}=\int_{t_1}^t\int_{\1 B_R}\frac{\1 G_R}{\1 n_y}(x,t,y,s)
u(y,s)\,d\sigma (y)\,ds\quad\forall 0<|x|<R, t_1<t\le t_2.\tag 4.3
$$
By (0.7) and Corollary 2.10, $\forall 0<\3<|x|<R$.
$$
|I_{2,\3}|
\le C\int_{t_1}^t\int_{\1 B_{\3}}
\frac{e^{-c|x-y|^2/(t-s)}}{(t-s)^{\frac{m+1}{2}}}|y|^{2-m}\,d\sigma (y)\,ds
\le\frac{C\3}{(|x|-\3)^{m+1}}.\tag 4.4
$$
Letting $\3\to 0$ in (4.1), by (4.2), (4.3), and (4.4),
$\forall 0<\3<|x|<R$, $t_1<t\le t_2$,
$$
u(x,t)
=\int_{B_R\setminus\{0\}}G_R(x,t,y,t_1)u(y,t_1)\,dy
+\int_{t_1}^t\int_{\1 B_R}\frac{\1 G_R}{\1 n_y}(x,t,y,s)
u(y,s)\,d\sigma (y)\,ds.\tag 4.5
$$
Since the right hand side of (4.5) is a $C^{\infty}$ function of $(x,t)\in
B_R\times (t_1,t_2)$, we can extend $u(x,t)$ to a continuous function
on $B_R\times (t_1,t_2)$ by defining $u(0,t)$ to be equal to the 
right hand side of (4.5). Then $u\in C^{\infty}(B_R\times (t_1,t_2))$
satisfies 
$$
u(x,t)
=\int_{B_R}G_R(x,t,y,t_1)u(y,t_1)\,dy
+\int_{t_1}^t\int_{\1 B_R}\frac{\1 G_R}{\1 n_y}(x,t,y,s)
u(y,s)\,d\sigma (y)\,ds.\tag 4.6
$$
Since $G_R(x,t,y,s)$ satisfies the heat equation, by (4.6) the 
extended $u$ satisfies the heat equation in $B_R\times 
(t_1,t_2)$. Since $0<t_1<t_2<T$ is arbitrary, $u$ has removable singularities
at $\{0\}\times (0,T)$ and the theorem follows.
\enddemo


\Refs

\ref
\key A\by D.G.~Aronson\paper Non-negative solutions of linear parabolic
equations\jour Ann. Scuola Norm. Sup. Pisa\yr 1968\vol 22(3)
\pages 607--694\endref

\ref
\key C\by B.~Chow\paper Lecture notes on Ricci flow I, II, III, 
Clay Mathematics Institute, Summer School\linebreak Program
2005 on Ricci Flow, 3-Manifolds and Geometry
June 20--July 16 at MSRI,\linebreak
http://www.claymath.org/programs/summer\_school/2005/program.php\#ricci
\endref

\ref
\key CK\by \ B.~Chow and D.~Knopf\book The Ricci flow:An introduction,
Mathematical Surveys and Monographs, Volume 110, Amer. Math. Soc.
\publaddr Providence, R.I., U.S.A.\yr 2004\endref

\ref
\key D\by D.M.~Deturck\paper Deforming metrics in the direction of their
Ricci tensors (improved version), in Collected Papers on Ricci Flow, 
ed. H.D.~Cao, B.~Chow, S.C.~Chu and S.T.~Yau, International Press, 
Somerville, MA, 2003\endref

\ref
\key DK\by \ B.E.J.~Dahlberg and C.~Kenig\paper Non-negative
solutions of generalized porous medium equations\jour Revista Mat.
Iberoamericana\vol 2\pages 267--305\yr 1986\endref

\ref
\key F\by A.~Friedman\book Partial Differential Equations of
Parabolic type\publ Prentice Hall, Inc.
\publaddr Englewood Cliffs, N.J.\yr 1964\endref

\ref
\key H1\by R.S.~Hamilton\paper Three-manifolds with positive Ricci curvature
\jour J. Differential Geom.\vol 17(2)\yr 1982\pages 255--306\endref

\ref
\key H2\by R.S.~Hamilton\paper Four-manifolds with positive curvature
operator\jour J. Differential Geom.\vol 24(2)\yr 1986\pages 153--179\endref

\ref
\key H3\by R.S.~Hamilton\paper The Ricci flow on surfaces\jour 
Contemp. Math.\vol 71\yr 1988\pages 237--261
\endref

\ref
\key H4\by R.S.~Hamilton\paper The Harnack estimate for the Ricci flow
\jour J. Differential Geom.\vol 37(1)\yr 1993\pages 225--243\endref

\ref 
\key H5\by R.S.~Hamilton\paper The formation of singularities in the Ricci flow
\jour Surveys in differential geometry, Vol. II (Cambridge, MA, 1993),7--136,
International Press, Cambridge, MA, 1995\endref

\ref
\key H6\by R.S.~Hamilton\paper A compactness property for solutions of the 
Ricci flow\jour Amer. J. Math.\vol 117(3)\yr 1995\pages 545--572\endref

\ref
\key Hs1\by \ \ S.Y.~Hsu\paper Global existence and uniqueness
of solutions of the Ricci flow equation\jour Differential
and Integral Equations\vol 14(3)\yr 2001\pages 305--320\endref

\ref
\key Hs2\by \ \ S.Y.~Hsu\paper Large time behaviour of solutions
of the Ricci flow equation on $R^2$\vol 197(1)\yr 2001
\pages 25--41\jour Pacific J. Math.\endref

\ref
\key Hs3\by \ \  S.Y.~Hsu\paper Asymptotic profile of
solutions of a singular diffusion equation as $t\to\infty$
\jour Nonlinear Analysis, TMA\vol 48\yr 2002\pages 781--790
\endref

\ref
\key Hs4\by \ \ S.Y.~Hsu\paper Dynamics of solutions of a
singular diffusion equation\jour Advances in Differential
Equations\vol 7(1)\yr 2002\pages 77--97\endref

\ref
\key Hs5\by \ \ S.Y.~Hsu\paper A simple proof on the non-existence of 
shrinking breathers for the Ricci flow\jour Calculus of Variations
and P.D.E.\vol 27(1)\yr 2006\pages 59--73\endref 

\ref
\key Hs6\by \ \ S.Y.~Hsu\paper Generalized $\Cal{L}$-geodesic and 
monotonicity of the generalized reduced volume in the Ricci flow,
http://arXiv.org/abs/math.DG/0608197
\endref 

\ref
\key Hs7\by \ \ S.Y.~Hsu, A lower bound for the scalar curvature of 
the standard solution of the Ricci flow, 
http://arXiv.org/abs/math.DG/0612426
\endref

\ref
\key Hu\by \ \ K.M.~Hui\paper A Fatou Theorem for the solution of the
heat equation at the corner points of a cylinder\jour Trans. A.M.S.
\vol 333(2)\yr 1992\pages 607--642
\endref

\ref
\key KL\by \ B.~Kleiner and J.~Lott\paper Notes on Perelman's papers,
http://arXiv.org/abs/math.DG/0605667\endref

\ref
\key LSU\by \ \ O.A.~Ladyzenskaya, V.A.~Solonnikov, and
N.N.~Uraltceva\book Linear and quasilinear equations of
parabolic type\publ Transl. Math. Mono. Vol 23,
Amer. Math. Soc.\publaddr Providence, R.I.\yr 1968\endref

\ref
\key LiT\by \ P.~Li and L.F.~Tam\paper The heat equation and harmonic maps
of complete manifolds\jour Invent. Math.\vol 105\yr 1991\pages 1--46
\endref

\ref
\key LT\by P. Lu and G.~Tian\paper Uniqueness of solutions in the work of
Perelman, http://www.math.lsa.umich.edu\linebreak
/$\sim$lott/ricciflow/StanUniqWork2.pdf
\endref

\ref
\key MT\by \ \ J.W.~Morgan and G.~Tian\paper Ricci flow and the Poincar\'e
Conjecture, http://arXiv.org/abs/math.DG\linebreak /0607607\endref

\ref
\key P1\by G.~Perelman\paper The entropy formula for the Ricci flow and its 
geometric applications,  http://arXiv.org\linebreak /abs/math.DG/0211159
\endref 

\ref
\key P2\by G.~Perelman\paper Ricci flow with surgery on three-manifolds,
http://arXiv.org/abs/math.DG/0303109\endref

\ref
\key S1\by W.X.~Shi\paper Deforming the metric on complete Riemannian manifolds
\jour J. Differential Geom.\vol 30\yr 1989\pages 223--301\endref

\ref
\key S2\by W.X.~Shi\paper Ricci deformation of the metric on complete 
non-compact Riemannian manifolds \jour J. Differential Geom.\vol 30\yr 1989
\pages 303--394\endref

\ref
\key W1\by \ L.F.~Wu\paper The Ricci flow on complete $R^2$
\jour Comm. in Analysis and Geometry\vol 1\yr 1993
\pages 439--472\endref

\ref
\key W2\by \ L.F.~Wu\paper A new result for the porous
medium equation\jour Bull. Amer. Math. Soc.\vol 28\yr 1993
\pages 90--94\endref

\ref
\key Ye\by Rugang Ye\paper On the $l$-function and the reduced volume of 
Perelman (preprint, Feb.,2004)\endref

\endRefs
\enddocument